\documentclass[12pt,a4paper]{article}
\usepackage{amsmath,amsxtra,amsfonts,amssymb,latexsym,amsthm,amscd,verbatim,amstext}
\usepackage[mathscr]{eucal}

\input cyracc.def

\newfont{\tricyr}{wncyr10 at 12pt}
\newfont{\tricyi}{wncyi10 at 12pt}
\newfont{\tricyb}{wncyb10 at 12pt}
\newfont{\Tricyr}{wncyr10 at 13.6pt}
\newfont{\Tricyi}{wncyi10 at 13.6pt}
\newfont{\Tricyb}{wncyb10 at 13.6pt}
\newfont{\tricmr}{cmr10 at 13.6pt}
\newfont{\tricmi}{cmti10 at 13.6pt}
\newfont{\tricmb}{cmb10 at 13.6pt}

\theoremstyle{plain}
\newtheorem{Th}{Theorem}
\newtheorem{Md}{Proposition}
\newtheorem{Lem}{Lemma}
\newtheorem{Cor}{Corollary}
\theoremstyle{definition}

\def\RR{\mathbb R}

\begin{document}
%%\NoBlackBoxes%%

\centerline {{\bf  GENERALIZED VARIATIONAL INEQUALITIES}} 
\centerline {{\bf  FOR MAXIMAL MONOTONE OPERATORS}}
\vskip 0.7cm
\begin{center}
Nga Quynh Nguyen$^{a,b}$
\end{center}
\begin{center}
{\footnotesize
$^a$ Graduate University of Science and Technology,\\
Vietnam Academy of Science and Technology,\\
 18 Hoang Quoc Viet, Cau Giay, Hanoi, Vietnam.\\
$^b$ Institute of Mathematics, 
Vietnam Academy of Science and Technology,\\
18 Hoang Quoc Viet, Cau Giay, Hanoi, Vietnam.}
\end{center}
\vskip 0.5cm 
{\bf Abstract}: In this paper we present some new results on the existence of solutions of generalized variational inequalities in  real reflexive Banach spaces with  Fr\'{e}chet differentiable norms. Moreover, we also  give some theorems about the structure of solution sets. The results obtained in this
paper improve and extend the ones announced by Fang and Peterson \cite{Fang} to infinite dimensional spaces.


\footnotetext[3]{{\it e-mail address}: nqnga@math.ac.vn}
\vskip 1.5cm 

\centerline{1. INTRODUCTION} 
\vskip 0.5cm 
\indent  Variational inequalities were initially introduced to deal with partial differential equations stimulated from mechanics (see, e.g., \cite{Hartman}).  They have been applied intensively to different fields such as mechanics, game theory, optimization, economics, etc.(see, e.g., \cite{Fang}-\cite{Browder1} and references therein). In \cite{Fang}, some results on the existence of solutions for generalized variational inequalities in ${\RR}^n$ for monotone operators, maximal monotone operators were established by Fang and Peterson. The results obtained in this paper improve and extend  their results to real reflexive Banach spaces with Fr\'{e}chet differentiable norms.  Moreover, we also give some theorems about the structure of solution sets. 
\vskip 0.5cm 
\centerline{2. PRELIMINARIES}
\vskip 0.5cm 

\indent Let $X$ be a real Banach space with the dual space $X^*$, $K$ be a subset in $X$. We denote by $\langle w^*, u \rangle$ the dual pair between  $w^*\in X^*$ and $u \in X$. Let $ T: K \rightarrow 2^{X^*} $ be a set-valued mapping. The generalized variational 
inequality problem, denoted by $GVI(K, T)$, is to find vectors  $\overline x\in K$, and $x^* \in T(\overline x)$ such that  
\[
\left\langle {\left. {x^* ,y - \overline x } \right\rangle } \right. \ge 0,\,\forall \,y \in K.
\]

Then $\overline x$ is called a solution of $GVI(K, T)$.  
Denote by $SOL(K, T)$ the solution set of $GVI(K, T)$. That is,

$$SOL(K,T)= \Bigl \{ \overline x\in K : \exists x^*\in T(\overline x)\; \text {such that} \left\langle {\left. {x^* ,y - \overline x } \right\rangle } \right. \ge 0,\,\forall \,y \in K \Bigr \}.$$

First we recall some concepts used in this work.

A mapping $T$ is called upper-semicontinuous at $x^0\in X$ if for every open subset $N$ containing $T(x^0)$ there exists an open neighborhood $M$ of $x^0$ such that $T(x)\subset N $ for all $x\in M$. $T$  is called upper-semicontinuous on $K$ if it is upper-semicontinuous at every $x^0\in K$. 
                                                                                                                                                                                                                                                                                                                                                                                                                                                                                                                                                                                                                                                                                                                                                                                                                                                                                               
Graph $G(T)$ is a subset of $K \times X^*$ defined by 
$$ G(T)=\lbrace(u,w): w\in T(u), u\in K \rbrace .
$$
We say that $T \subseteq T_1$ if $G(T)\subseteq G(T_1)$. 
We say that
\begin{enumerate}
\item [(i)] $T$ is pseudomonotone on $K$ if from 
$\left\langle {\left. { x^* , \;  y - x } \right \rangle } \right.  \ge 0$, it follows that  
$\left\langle {\left. {y^* ,\;  y - x } \right \rangle } \right. \ge 0$
for every $x,y \in K, x^* \in T(x), y^* \in T(y)$;

\item [(ii)] $T$  is monotone on  $K$ if for every $x,y \in K$ and $ x^*\in T(x), y^* \in T(y)$, we have $\left\langle {\left. { y^* - x^* ,\; y -  x } \right \rangle } \right.  \ge 0$;

\item [(iii)] $T$  is  strictly monotone on $K$ if for every $x,y \in K, x\ne y$ and $x^*\in T(x), y^* \in T(y)$, we have  $\left\langle {\left. { y^* - x^* ,\; y -  x } \right \rangle } \right.  > 0$;

\item [(iv)] $T$ is strongly monotone on  $K$ if there exists a number $\alpha >0 $ such that for every $x,y \in K$ and $x^*\in T(x), y^* \in T(y)$, we have  \\ $\left\langle {\left. { y^* - x^* ,\; y -  x } \right \rangle } \right.  \ge \alpha || x-y ||^2  $;

\item [(v)] $T$ is maximal monotone on $K$ if  $T$ is monotone on $K$ and for each monotone mapping $T_1$  such that $T\subseteq T_1$, we have $T=T_1$.
\end{enumerate}

From the above definition we note that strongly monotonicity $\Rightarrow$ strictly monotonicity
$\Rightarrow $ monotonicity $\Rightarrow$ pseudomonotonicity.

A cone in  $X$ is nonempty subset $K\subset X$ such that $tx\in K$ for all $t\ge 0$  and  $x\in K$. A cone $K$ is called pointed if $K\cap (-K)= \{0 \}$.
We denote $$ K^* =\lbrace \lambda ^* \in X^*: \langle  \lambda ^* ,\; x \rangle \ge 0 \; \text {for all}\; x\in K \rbrace . $$  
Obviously that $K^*$ is a convex, closed cone in $X^*$.

A subset is said to be solid if its interior is nonempty. The notations  $\text{int}(K), \partial K$ mean the interior and the boundary of $K$, respectively. A subset $K\subset X$ is said to be contractible if there exist a point $x^0 \in K$ and  a continuous function $g : K\times [0,1] \rightarrow K$ such that $g(x,0)=x$ and\; $g(x,1)=x^0$ for all $x\in K$. Note that if $K$ is convex then $K$ is contractible.

Now we recall some previous  results which we will use in this paper.

\begin{Lem}[see \cite{Junlouchai}]\label{l1}
Let $X, Y$ be two metric spaces and $T: X \rightarrow 2^Y$ be a set-valued mapping. Given any $x\in X$, if $T(x)$ is compact and $ T $ is upper-semicontinuous at $x$ then for every sequence $\{x_n\}$  converging to $x$, every sequence $u^{*}_{n}\in T(x_n)$ must have a cluster point $u^*\in T(x)$.
\end{Lem}

\begin{Md}[see \cite{Aubin}]\label{p1} Let $X, Y$ be two Hausdorff topological vector spaces and $T: X \rightarrow 2^Y $ be an upper-semicontinuous map with nonempty compact values. If $X$ is a compact space then $T(X)$ is compact.
\end{Md}

\begin{Lem}[see \cite{Browder2}]\label{l2}
Let $X$ be a real reflexive Banach space and $T: X \rightarrow 2^{X^*} $ be a maximal monotone map. Then  
\begin{enumerate}
\item [(i)] For each $u\in X, T(u)$ is a closed convex subset in $X^*$ ;
\item [(ii)] If $\left\{u_k\right\},\left\{ v_k\right\} $ are two sequences in $X$ and ${X^*}$, respectively, such that $u_k \rightarrow \overline u, v_k \in T(u_k)$ and $v_k \rightharpoonup \overline v$, then $\overline v\in T(\overline u)$.
\end{enumerate}
\end{Lem}

\begin{Lem}[see \cite{Browder2}]\label{l3}
Let $X$ be a real reflexive Banach space and $T: X \rightarrow 2^{X^*} $be a maximal monotone map. If $\left\{u_k\right\}$, $\left\{ v_k\right\} $ are two sequences in $X$ and ${X^*}$, respectively, such that $u_k \rightharpoonup \overline u$, $v_k \in T(u_k)$ and $v_k \rightarrow \overline v $, then $\overline v\in T(\overline u)$.
\end{Lem}

\begin{Th}[see \cite{Fang}]\label{t1}(Hartman-Stampachia, Saigal)
Suppose that
\begin{enumerate}
\item [(i)] $K$ is a nonempty, compact, convex subset in ${\RR}^n$;
\item [(ii)] $T: K \rightarrow 2^ {{\RR}^n}$ is an upper-semicontinuous;
\item [(iii)] $T(x)$ is a nonempty, compact and contractible set in ${\RR}^n$ for each $x \in K$.

Then $GVI(K, T)$ has a solution.
\end{enumerate}
\end{Th}

\begin{Th}[see \cite{Junlouchai}]\label{t2}
Let $X$  be a real reflexive Banach space with a Fr\'{e}chet differentiable norm. Suppose that
\begin{enumerate}
\item [(i)] $K$ is a nonempty, compact, convex subset in $X$
\item [(ii)] $T: K \rightarrow 2^{X^*} $ is an upper-semicontinuous;

\item [(iii)] $T(x)$ is a nonempty, closed and contractible subset in $X^*$ for each \\ $x \in K$;

\item [(iv)] $T(K)=\underset {x\in K}\bigcup T(x)$ is compact in $X^*$. 

Then $GVI(K, T)$  has a solution.
\end{enumerate}
\end{Th}

The following theorem gives the property of interior points of $K^*$.
\begin{Th}[see \cite{Chiang}]\label{t3}
Let $X$ be a locally convex Hausdorff topological vector space and $K $ be a pointed cone. Then $\lambda ^*\in $ {\rm int}$(K^*)$  if and only if $K(\lambda ^*,r)$ is bounded for all $r > 0$, where $K(\lambda ^*,r)=\lbrace x\in K: \langle \lambda ^*, x \rangle \le r \rbrace $. 
\end{Th}

\vskip 0.5cm 
\centerline{3. MAIN RESULTS} 
\vskip 0.5cm 
\begin{Th}\label{t4}
Let $X$ be a real reflexive Banach space and $K$ be a nonempty subset in $X$. For a given maximal monotone  map $ T: K \rightarrow 2^{X^*} $ on $K$, putting $$ S = \lbrace (x,y^*)\in K\times X^*: y^*\in T(x) \;  \text {and} \; \langle \;  y^*, z-x \rangle \ge 0, \; \forall z\in K \rbrace. $$ Then 
\begin{enumerate}
\item [(i)] if $K$  is convex then $S$ is also convex.

\item [(ii)]if $K$ is closed then $S$ is also closed.

\item [(iii)] if $K$ is compact and $T(K)$ is weakly compact then $S$ is compact.

\item [(iv)] if $K$ is weakly compact and $T(K)$ is compact then $S$ is compact.

\end{enumerate}

\end{Th}

{\it Proof.} $(i)$ If $S=\emptyset$ or $S$ contains only one point then it is obviously true. Therefore we can suppose that $S\ne \emptyset$ and $(x_1,y_1^*), (x_2,y_2^*)\in S$. Let $\lambda_1,\lambda_2\in [0,1]$ such that $\lambda_1+\lambda_2 = 1.$ Since $K$  is convex, $\lambda_1x_1+\lambda_2x_2\in K.$ For every $z\in K$ we have $\langle \;  y_1^*, z-x_1 \rangle \ge 0$  and  $ \langle \; y_2^*, z-x_2 \rangle \ge 0$. Hence by monotonicity of $T$, we have 
\begin{eqnarray*}
\langle \lambda_1 y_1^* + \lambda_2 y_2^* , z - (\lambda_1 x_1 + \lambda_2 x_2) \rangle &=& \langle \lambda_1 y_1^* + \lambda_2 y_2^*,\lambda_1(z - x_1) + \lambda_2(z- x_2)\rangle \\
&\ge &\lambda_1\lambda_2 \big ( \langle y_2^*, z - x_1  \rangle + \langle y_1^*, z - x_2 \rangle \big )\\
&=& \lambda_1\lambda_2 \big ( \langle y_2^*, z \rangle + \langle y_1^*, z \rangle - \langle y_2^*, x_1 \rangle -\langle y_1^*, x_2 \rangle  \big )\\
&\ge & \lambda_1\lambda_2 \big ( \langle y_2^*, x_2 \rangle + \langle y_1^*, x_1 \rangle - \langle y_2^*, x_1 \rangle -\langle y_1^*, x_2 \rangle  \big )\\
& = & \lambda_1\lambda_2 \langle y_1^*-y_2^*, x_1-x_2 \rangle \\
&\ge & 0.
\end{eqnarray*}
Now we prove that $\lambda_1 y_1^* + \lambda_2 y_2^*\in T(\lambda_1 x_1 + \lambda_2 x_2) $.\\  We consider the map $T_1$ defined by

\noindent $T_1(x):= T(x)\;\bigcup \; \big \lbrace \; y^*\in X^* : \exists  (x_1,y_1^*), (x_2,y_2^*)\in S $ \text {и} 
$\lambda_1, \lambda_2 \ge 0, \lambda_1+\lambda_2 = 1$
\indent \indent \indent \indent \indent \indent \indent \text {such that} $x=\lambda_1x_1+\lambda_2x_2,\; y^*=\lambda_1y_1^*+ \lambda_2y_2^*\big \rbrace $ \\
for each $x\in K$. Obviously, $T\subseteq T_1$. Now we prove $T_1$ is monotone. Let $x_1,x_2\in K$  and $y_1^*\in T_1(x_1), y_2^*\in T_1(x_2)$.

\noindent {\bf Case 1: } If $y_1^*\in T(x_1)$ and $y_2^*\in T(x_2)$, the monotonicity of $T$ implies that $\langle \; y_1^*-y_2^*\; , x_1-x_2 \rangle \ge 0.$

\noindent {\bf Case 2: }If $y_1^*\in T(x_1)$ and $x_2 = \lambda_1 u + \lambda_2 v,\;  y_2^* = \lambda_1 u^*+ \lambda_2 v^* $, for $(u,u^*), (v,v^*)\in S$ and $\lambda_1, \lambda_2 \ge 0, \lambda_1+\lambda_2 = 1$, from the monotonicity of $T$ and $y_1^*\in T(x_1), u^*\in T(u), v^*\in T(v)$, we have 
\allowdisplaybreaks
\begin{eqnarray*}
\langle \; y_2^*-y_1^*\; , x_2-x_1 \rangle & = & \langle \;  \lambda_1 u^*+ \lambda_2 v^* - y_1^*\; , \lambda_1 u+\lambda_2 v-x_1 \rangle \\
& = & \langle \;\lambda_1 (u^* - y_1^*) + \lambda_2(v^*-y_1^*) \; , \lambda_1 (u- x_1) + \lambda_2(v -x_1) \rangle  \\
& = & \lambda_1^2 \langle \; u^*-y_1^*\; , u - x_1  \rangle + \lambda_2^2 \langle \; v^*-y_1^*\; , v - x_1  \rangle \\
& + & \lambda_1 \lambda_2\big [\langle \; u^*-y_1^*\; , v - x_1  \rangle + \langle \; v^*-y_1^*\; , u - x_1  \rangle \big ] \\
& = & \lambda_1 \big [ \lambda_1  \langle \; u^*-y_1^*\; , u - x_1  \rangle + \lambda_2 \langle \; v^*-y_1^*\; , u - x_1  \rangle \big ] \\
& + & \lambda_2 \big [ \lambda_2  \langle \; v^*-y_1^*\; , v - x_1  \rangle + \lambda_1 \langle \; u^*-y_1^*\; , v - x_1  \rangle \big ] \\
& = & \lambda_1 \big [ (1 - \lambda_2)  \langle \; u^*-y_1^*\; , u - x_1  \rangle + \lambda_2 \langle \; v^*-y_1^*\; , u - x_1  \rangle \big ] \\
& + & \lambda_2 \big [ (1- \lambda_1)  \langle \; v^*-y_1^*\; , v - x_1  \rangle + \lambda_1 \langle \; u^*-y_1^*\; , v - x_1  \rangle \big ] \\
& = & \lambda_1 \big [ \langle \; u^*-y_1^*\; , u - x_1  \rangle - \lambda_2 \langle \; u^*- v^*\; , u - x_1  \rangle \big ] \\
& + & \lambda_2 \big [ \langle \; v^*-y_1^*\; , v - x_1  \rangle - \lambda_1 \langle \; v^*- u^*\; , v - x_1  \rangle \big ] \\
& = & \lambda_1 \langle \; u^*-y_1^*\; , u - x_1  \rangle + \lambda_2  \langle \; v^*-y_1^*\; , v - x_1  \rangle \\
& - & \lambda_1  \lambda_2 \big [ \langle \; u^*- v^*\; , u - x_1  \rangle + \langle \; v^*- u^*\; , v - x_1  \rangle \big ] \\
&\ge & -\lambda_1  \lambda_2  \langle \; u^*- v^*\; , u - v  \rangle.
\end{eqnarray*}
Since $(u,u^*), (v,v^*)\in S$, we have $ \langle \; u^*\; , v - u  \rangle \ge  0$ and $ \langle \; v^*\; , u - v  \rangle \ge  0.$ Therefore, $\langle \; u^*- v^*\; , u - v  \rangle \le 0$. On the other hand, since $T$ is monotone, $\langle \; u^*- v^*\; , u - v  \rangle \ge 0.$ Thus, $\langle \; u^*- v^*\; , u - v  \rangle = 0$ and hence $\langle \; y_2^*-y_1^*\; , x_2-x_1 \rangle \ge 0 $.

\noindent {\bf Case 3:} If $y_2^* \in T(x_2)$ and $x_1 = \lambda_1 u + \lambda_2 v,\;  y_1^* = \lambda_1 u^*+ \lambda_2 v^* $ for $(u,u^*), (v,v^*)\\\in S $ and $\lambda_1, \lambda_2 \ge 0, \lambda_1+\lambda_2 = 1 $, then by similar arguments as in Case 2 we also obtain $\langle \; y_2^*-y_1^*\; , x_2-x_1 \rangle \ge 0 $.

\noindent {\bf Case 4:} If $x_1 = \lambda_1 u_1 + \lambda_2 v_1,\;  y_1^* = \lambda_1 u_1^*+ \lambda_2 v_1^* $ for $(u_1,u_1^*), (v_1,v_1^*)\in S $ and $\lambda_1, \lambda_2 \ge 0, \lambda_1+\lambda_2 = 1 $, and $x_2 = \mu_1 u_2 + \mu_2 v_2,\;  y_2^* = \mu_1 u_2^*+ \mu_2 v_2^* $ for $(u_2,u_2^*), (v_2,v_2^*)\in S $ and $\mu_1, \mu_2 \ge 0, \mu_1+\mu_2 = 1 $, then
\begin{eqnarray*}
\langle y_2^*-y_1^*,x_2-x_1\rangle &=& \langle \mu_1 u_2^*+ \mu_2 v_2^* - \lambda_1 u_1^*-\lambda_2 v_1^*,\mu_1 u_2 +\mu_2 v_2-\lambda_1 u_1 -\lambda_2v_1 \rangle \\
&=&\mu_1^2 \langle u_2^*, u_2 \rangle + \mu_2^2 \langle v_2^*, v_2 \rangle + \lambda_1^2 \langle u_1^*, u_1\rangle + \lambda_2^2 \langle v_1^*,v_1 \rangle \\
&+&\mu_1\mu_2 \big [\langle u_2^*, v_2 \rangle + \langle v_2^*,u_2 \rangle \big ]-\lambda_1\mu_1 \big [\langle u_1^*, u_2 \rangle + \langle u_2^*, u_1 \rangle \big ]\\
&-&\lambda_1\mu_2 \big [\langle u_1^*, v_2 \rangle + \langle v_2^*,u_1 \rangle \big ]-\lambda_2\mu_1 \big [\langle u_2^*, v_1 \rangle + \langle v_1^*,u_2 \rangle \big ]\\
& - & \lambda_2\mu_2 \big [ \langle v_2^*, v_1 \rangle + \langle v_1^*, v_2 \rangle \big ] +  \lambda_1\lambda_2 \big [\langle u_1^*, v_1 \rangle + \langle v_1^*, u_1 \rangle \big ]\\
&=&\mu_1^2\langle u_2^*,u_2\rangle + \mu_2^2 \langle v_2^*, v_2 \rangle + \lambda_1^2 \langle u_1^*, u_1 \rangle + \lambda_2^2 \langle v_1^*, v_1 \rangle\\
& + & \mu_1\mu_2 \big [\langle u_2^*, u_2 \rangle + \langle v_2^*, v_2 \rangle \big ]-\lambda_1\mu_1 \big [\langle u_1^*, u_1 \rangle + \langle u_2^*, u_2 \rangle \big ]\\
& - & \lambda_1\mu_2 \big [\langle u_1^*, u_1 \rangle + \langle v_2^*, v_2 \rangle \big ] - \lambda_2\mu_1 \big [\langle u_2^*, u_2 \rangle + \langle v_1^*, v_1 \rangle \big ]\\
& - & \lambda_2\mu_2 \big [ \langle v_2^*, v_2 \rangle + \langle v_1^*, v_1 \rangle \big ] +  \lambda_1\lambda_2 \big [\langle u_1^*, u_1 \rangle + \langle v_1^*, v_1 \rangle \big ]\\
&=&\big(\lambda_1^2-\lambda_1\mu_1-\lambda_1\mu_2+\lambda_1 \lambda_2 \big) \langle u_1^*,u_1\rangle \\
&+& \big(\mu_1^2+ \mu_1\mu_2- \lambda_1\mu_1 -  \lambda_2\mu_1\big)\langle u_2^*,u_2\rangle \\
&+ &\big (  \mu_2^2 + \mu_1\mu_2- \lambda_1\mu_2 -  \lambda_2 \mu_2 \big ) \langle v_2^*, v_2 \rangle\\
&+& \big (\lambda_2^2- \lambda_2\mu_1- \lambda_2\mu_2 + \lambda_1 \lambda_2 \big ) \langle v_1^*, v_1 \rangle\\
&=& \big [ \lambda_1 \big (\lambda_1 + \lambda_2 \big ) - \lambda_1 \big ( \mu_1 + \mu_2 \big )  \big ] \langle u_1^*, u_1 \rangle \\
&+& \big [ \mu_1 \big (\mu_1 + \mu_2 \big ) - \mu_1 \big ( \lambda_1 + \lambda_2 \big )  \big ] \langle u_2^*, u_2 \rangle \\
&+& \big [ \mu_2 \big (\mu_1 + \mu_2 \big ) - \mu_2 \big ( \lambda_1 + \lambda_2 \big )  \big ] \langle v_2^*, v_2 \rangle \\
&+& \big [ \lambda_2 \big (\lambda_1 + \lambda_2 \big ) - \lambda_2 \big ( \mu_1 + \mu_2 \big )  \big ] \langle v_1^*, v_1 \rangle \\
&=& 0.
\end{eqnarray*}
Therefore $T_1$ is monotone. The maximal monotonicity of $T$ implies that $T=T_1$ and $\lambda_1 y_1^* + \lambda_2 y_2^* \in T(\lambda_1 x_1+ \lambda_2 x_2).$  Hence $S$ is convex.

$(ii)$ In this part we prove that $S$ is closed. Let $(x_n,y_n^*)\in S $ and $(x_n,y_n^*)\rightarrow (x,y)$. Since $(x_n,y_n^*)\in S$, we have $x_n\in K, y_n^*\in T(x_n)$ and $\langle y_n^*\;,z-x_n\rangle \ge 0 $ for all $z\in K$. By the maximal monotonicity of $T$ and $x_n \rightarrow x,y_n^* \rightarrow y^*, y_n^*\in T(x_n)$ and using Lemma \ref{l2}, we imply $y^*\in T(x)$. Because $\langle y_n^*\;,z-x_n\rangle \ge 0 $  for all $z\in K$ and $x_n \rightarrow x,y_n^* \rightarrow y^*$ we obtain $\langle y^*\;,z-x\rangle \ge 0 $ for  $z\in K$. Since $K$ is closed, $x\in K$. Thus, $(x,y^*)\in S$ and hence $S$ is closed.

$(iii)$ Let $(x_n,y_n^*)$ be an arbitrary sequence in $S$. Then $x_n\in K$ and $y_n^*\in T(x_n)\subset T(K)$ and $\langle y_n^*\;,z-x_n\rangle \ge 0 $ for all $z\in K$. Since $K$ is compact and $T(K)$ is weakly compact, there exists a subsequence $(x_n,y_n^*)$ such that $x_n$ converges to $x\in K$, and $y_n^*$ weakly converges to $y^*\in T(K)$. By using similar arguments as in $(ii)$ and using Lemma \ref{l2} we obtain the desired conclusion.

$(iv)$ By using similar arguments as in $(iii)$ and using Lemma \ref{l3} we obtain the desired conclusion.
\qed

From now on we assume that $X$ is a real reflexive Banach space with a Fr\'{e}chet differentiable norm.

The following lemma extends the theorem of Hartman- Stampachia, Saigal for real reflexive Banach spaces with Fr\'{e}chet differentiable norms.

\begin{Lem}\label{l4}
Suppose that 
\begin{enumerate}
\item [(i)] $K$ is a nonempty compact convex subset in $X$;

\item [(ii)] $T: K \rightarrow 2^{X^*}$ is upper-semicontinuous;

\item [(iii)] $T(x)$ is a nonempty compact contractible subset in $X^*$ for each $x \in K$.

Then $SOL(K, T)$ is nonempty and compact.
\end{enumerate}
\end{Lem}

{\it Proof.} By Proposition \ref{p1}, $T(K)$ is compact. By Theorem \ref{t2}, $SOL(K,T)$ is nonempty. Assume that $\left\{x_n\right\}$ is a sequence in $SOL(K,T)$ and $x_n \rightarrow x.$ Since $\left\{x_n\right\}\in $ 
$SOL(K,T)$ there exist $u_n^*\in T(x_n)$ such that $\langle u_n^*\;, y-x_n\rangle \ge 0$ for all $y\in K$. Because $T(x)$ is compact and 
$T$ is upper-semicontinuous at $x$, by Lemma \ref{l1}, the sequence $\left\{u_n^*\right\}$ has a cluster point $u^*\in T(x)$. Without loss of generality, we can assume that $u_n^*\rightarrow u^*\in T(x)$. From $\langle u_n^*\;, y-x_n\rangle \ge 0$ for all $y\in K$, it follows that $\langle u^*\;, y-x\rangle \ge 0$ for all $y\in K$. Thus $x\in $ $SOL(K,T)$ and hence $SOL(K,T)$ is closed. Because $SOL(K,T)\subset K$ and $K$ is compact we can conclude that
$SOL(K,T)$ is compact.\qed

We note that Lemma \ref{l4} has stronger conclusion than the conclusion of Theorem \ref{t2}. Namely, we can conclude about the structure of solution set.

\begin{Lem}\label{l5}
Let $K$ be a nonempty convex subset in $X$ and $C$ be a solid convex subset in $X$. Assume that there exist $x^0\in K\cap $ {\rm int}$(C)$ and $y^*\in X^*$ such that $\langle y^*\;, x-x^0\rangle\ge 0$ for all $x\in K\cap C$. Then $\langle y^*\;, x-x^0\rangle\ge 0$ for all $x\in K$.
\end{Lem}

{\it Proof.}
For each $x \in K$ there exists a small enough number $s\ge 0$ such that $z=(1-s)x^0 + sx\in K\cap C$. Thus, $\langle y^*, z-x^0\rangle\ge 0$. Therefore, $\langle y^*, x-x^0\rangle\ge 0$.\qed

\begin{Th}\label{t5}
Let $K $ be a nonempty convex subset in $X$  and $ T: K \rightarrow 2^{X^*}.$  Suppose that there exists a solid  convex subset $C$ in $X$ such that

\begin{enumerate}
\item [(i)] $K \cap C $ is nonempty and compact;
\item [(ii)] $T\Big |_{K\cap C}$ is upper-semicontinuous;
\item [(iii)] $T(x)$ is nonempty compact contractible subset for all  $x\in K\cap C$;
\item[(iv)] For each $x\in K\cap \partial C$, there exists $x^0 \in K\cap {\rm int}(C)$ such that
  $0 \le \langle y^*, x - x^0 \rangle  $ for all  $y^*\in T(x)$.
\end{enumerate}
Then $SOL(K,T)\cap C$ is nonempty and compact.
\end{Th}

{\it Proof.}
By Lemma \ref{l4}, $SOL(K\cap C,T)$ is nonempty.  Let $\overline x \in SOL(K\cap C,T)$. Then $\overline x \in K\cap C$ and there exists  $y^*\in T(\overline x)$ such that 
\begin{equation}\label{eq:1}
\langle y^*, x - \overline x \rangle \ge 0 
\end{equation}
for all $x\in K\cap C$.

From $\overline x\in K\cap C$ it follows that either  $\overline x\in K\cap \partial C$ or $\overline x\in K\cap {\rm int}(C)$. If $\overline x\in K\cap {\rm int}(C)$, then putting $x_0 = \overline x$ in Lemma \ref{l5}, we can conclude that $\langle y^*, x - \overline x \rangle \ge 0$  for all $x\in K$. Therefore,
$\overline x\in SOL(K,T)\cap C$ and hence $SOL(K,T)\cap C$ is nonempty. If $\overline x\in K\cap \partial C$, then by assumption $(iv)$ there exists $x^0 \in K\cap {\rm int}(C)$ such that
\begin{equation}\label{eq:2}
\langle y^*, \overline x -x^0 \rangle \ge 0. 
\end{equation}
From \eqref{eq:1} and \eqref{eq:2} we imply that $\langle y^*, x -x^0 \rangle \ge 0$ for all $x\in K\cap C$. Again by Lemma \ref{l5}, we obtain 
\begin{equation}\label{eq:3}
\langle y^*, x - x^0 \rangle \ge 0 
\end{equation}
for all $x\in K$.
Because $x^0 \in K\cap C $, from \eqref{eq:1} it follows that
\begin{equation}\label{eq:4}
\langle y^*, x^0 - \overline x \rangle \ge 0. 
\end{equation}
From \eqref{eq:3} and \eqref{eq:4} we get $\langle y^*, x - \overline x \rangle \ge 0$ for all $x\in K$. Thus,
$\overline x\in SOL(K,T)\cap C$ and hence $SOL(K,T)\cap C$ is nonempty. 

Let $x_n \in SOL(K,T) \cap \; C$  and $x_n \rightarrow x$. Because $x_n \in K \cap \; C$  and $K \cap \; C$ is compact, we have $x \in K\cap \; C $. Since $x_n \in SOL(K,T)$,  there exists $u_n^* \in T(x_n)$ such that $\langle u_n^* , y-x_n \rangle \ge 0$ for all $y\in K$. Because $T(x)$ is compact and $T$ is upper-semicontinuous at $x$, by Lemma \ref{l1}, the sequence $\left\{u_n^*\right\}$ has a cluster point $u^*\in T(x)$. Without loss of generality, we can assume that $u_n^* \rightarrow u^* \in T(x)$. Therefore $\langle u^*, y - x \rangle \ge 0$ for all $y\in K$. Thus, $x\in SOL(K,T)$ and hence $SOL(K,T)\cap C$ is closed. Since $SOL(K,T)\cap C\subset K\cap \;C$ and  $K\cap C$ is compact, we can conclude that $SOL(K,T)\cap C$ is compact.\qed

Put

$$    B_r = \Bigl \{  x\in X: || x||  \le r    \Bigr \}, \hspace {.3cm }    C_r = \Bigl \{  x\in X: || x||  = r    \Bigr \}.   $$

For arbitrary nonempty subset $K$ in $X$, put 
$$    B(K)  = \Bigl \{  r\ge 0:  K\cap B_r \ne\emptyset    \Bigr \}, \hspace {.3cm }    C(K) = \Bigl \{  r\ge 0:  K\cap C_r \ne\emptyset    \Bigr \}.   $$

It is obvious that $ C(K)\subset B(K) $ and if $K$ is nonempty convex then $C(K)$ is a nonempty interval. If $K$ is not bounded then $C(K)$ is also not bounded.

For a given nonempty convex subset $K$ in $X$ and a set-valued mapping with nonempty values $T: K \rightarrow 2^{X^*} $, we define a function $ C_T(r, x^0)$ as follows:
$$C_T: \bigl ( B(K) \backslash \{0\} \bigr ) \times X \longrightarrow \RR \cup \{-\infty  \},$$
\[  C_T(r, x^0) = \begin{cases}
\underset {x\in K \cap C_r}{\inf}\Bigl(\underset {z^*\in T(x)}{\inf} {\langle z^*, x - x^0\rangle}/\bigl ({r +  ||x^0 ||}\bigr)\Bigr), & \, \text {if}\; r\in C(K)\\
0, &\, \text {if}\; r\in B(K)\backslash  C(K). 
\end{cases} \]

\begin{Th}\label{t6}
Suppose that
\begin{enumerate}
\item [(i)] $K$ is a nonempty convex closed subset in $X$ such that every weakly convergent sequence in $K$ strongly convergent, $T: K \rightarrow 2^{X^*}$ is upper-semicontinuous, where $T(x)$ is a nonempty compact contractible subset for every  $x\in K$;
\item [(ii)] There exist $r > 0$ and $x^0 \in K \cap $ {\rm int}$(B_r)$ such that $0 \le C_T(r, x^0).$
\end{enumerate}
Then, for each given element $q^* \in X^*$ satisfying $||q^* || \le C_T(r, x^0), $ the solution set $SOL(K, T + q^*)\cap B_r$ is nonempty and compact, where 
\[  ( T + q^* ) (x) = \bigl\{   z^* + q^* : \, z^* \in T(x)                \bigr \}.   \] 
\end{Th}

{\it Proof.}
Because $x^0 \in K \cap {\rm int}(B_r)$ we have $x^0 \in K \cap B_r$. Thus, $K \cap B_r$ is nonempty. Since $K \cap B_r$ is convex closed, it must be weakly closed. Since $X$ is a reflexive Banach space, it follows that $B_r$ is weakly compact. Therefore, $K \cap B_r$ is weakly compact. By assumption $(i)$  every weakly convergent sequence in $K$ strongly convergent, it follows that $K \cap B_r$  is compact. Let $x \in K\cap C_r$ and $y^*=z^*+q^*$, where $z^*\in T(x)$. Then
\begin{eqnarray*}
\langle y^*, x - x^0 \rangle & = & \langle z^*, x - x^0 \rangle + \langle q^*, x - x^0 \rangle \\
& \ge & (r +  ||x^0 ||)C_T(r, x^0)- ||q^* ||||x-x^0 ||\\
& \ge & (r +  ||x^0 ||)C_T(r, x^0)- (||x||+||x^0 ||)||q^* ||\\
& = & (r +  ||x^0 ||)(C_T(r, x^0)- ||q^* ||)\\
& \ge & 0.
\end{eqnarray*}
By Theorem \ref{t5}, $SOL(K,T+ q^*)\cap B_r$ is nonempty and compact.\qed

\begin{Th}\label{t7}
Suppose that the condition $(i)$ in Theorem \ref{t6} is satisfied. Moreover, suppose that
\begin{enumerate}
\item [(ii')] There exist $r > 0$ and  $x^0\in K\cap$  {\rm int}$(B_r)$ such that           
\[   0 < \beta := \underset {r' \ge r} {\inf} C_T(r', x^0) .             \]
\end{enumerate}
Then for each given element $q^* \in X^*,$ satisfying $||q^*|| < \beta, $ the solution set SOL$(K, T + q^*)\cap B_r$ is nonempty and compact. Moreover, $SOL(K,T+ q^*)\subset {\rm int}(B_r)$.
\end{Th}

{\it Proof.}
By Theorem \ref{t6}, $SOL(K,T+ q^*)\cap B_r$ is nonempty and compact. Now we prove that $SOL(K,T+ q^*)\subset {\rm int}(B_r)$.
Suppose that $\overline x $ is a solution of $GVI(K,T+q^*)$. Set $\overline r = ||\overline x||$. If $\overline r = 0$, then it is obvious that $||\overline x||=0 < r$. If $\overline r > 0$, then since $\overline x \in SOL(K,T+q^*)$, there exists $y^*=z^*+q^*$ where $z^*\in T(\overline x)$ such that $\langle y^*, x - \overline x \rangle \ge 0$ for all $x\in K$ or say in other words, $\langle z^*+q^* , x - \overline x \rangle \ge 0$  for all  $x\in K$. Because $x_0\in K$ we get $\langle z^*+q^* , x^0 - \overline x \rangle \ge 0$. Thus, 
$\langle z^* , \overline x - x^0  \rangle \le \langle q^* , x^0- \overline x  \rangle \le ||x^0-\overline x|| ||q^*|| \le  (||x^0|| + ||\overline x||)||q^*|| = (||x^0||+ \overline r)||q^*||$. Hence, $C_T(\overline r, x^0)\le ||q^*||$. Since $||q^*||<\beta =  \underset {r' \ge r} {\inf} C_T(r', x^0)$ we have $\overline r < r$. Therefore, $||\overline x|| < r$.\qed

\begin{Cor}\label{c1}
{\it Suppose that the condition $(i)$ in Theorem \ref{t6} is satisfied. Moreover, suppose that}
\begin{enumerate}
\item [($ii^{''}$)] {\it There exists $x^0 \in K, $ such that}
\[
\mathop {\lim }\limits_{\scriptstyle || x || \to  + \infty  \hfill \atop 
  \scriptstyle x \in K \hfill} \bigl[ \underset{z^*\in T(x)} {\inf}{\langle z^*, x - x^0 \rangle}/ {||x||}    \bigr] = + \infty.\]
{\it Then for an arbitrary element  $q^* \in X^*$, the solution set SOL$(K, T + q^*)$ is nonempty and compact. Moreover, there exists $r>0$ such that \\$SOL(K,T+ q^*)\subset {\rm int}(B_r)$ for each element $q^* \in X^*$.}
\end{enumerate}
\end{Cor}

{\it Proof.}
Choose $r$ such that $r> ||x^0||$  and when $||x||\ge r$

$$\underset{z^*\in T(x)} {\inf}{\langle z^*, x - x^0 \rangle}/ {||x||}  > 2(||q^*||+1).$$

Since $r> ||x^0||$ and $||x||\ge r$ we have 
$||x|| > ||x^0||$. Because $r> ||x^0||$, it follows that $x^0\in K\cap \; {\rm int}(B_r).$ For $r'\ge r$  we have
\begin{eqnarray*}
C_T(r', x^0)&=&\underset {x\in K \cap C_{r'}} {\inf} \;\;      \underset {z^*\in T(x)} {\inf} \frac{\langle z^*, x - x^0\rangle}{r' +  ||x^0 ||}\\
& = &\underset {x\in K \cap C_{r'}} {\inf} \;\;      \underset {z^*\in T(x)} {\inf} \frac{\langle z^*, x - x^0\rangle}{||x|| +  ||x^0 ||}.
\end{eqnarray*}
Since
\begin{eqnarray*}
\frac{\langle z^*, x - x^0\rangle }{||x|| +  ||x^0 ||} & > &\frac{2(||q^*||+1)||x||}{||x|| +  ||x^0 ||}\\
& > & \frac{(||q^*||+1)(||x|| +  ||x^0 ||)}{||x|| +  ||x^0 ||}\\
&=& ||q^*||+1, 
\end{eqnarray*}
for all  $x\in K\cap C_{r'}$  and $z^*\in T(x)$  we have $C_T(r', x^0)\ge ||q^*||+1>0$. Thus,  it follows that $\underset {r' \ge r} {\inf} C_T(r', x^0)\ge ||q^*||+1>0$. Now the Corollary \ref{c1} follows from Theorem \ref{t7}.\qed

Note that the condition $(ii {}^{''})$ is the traditional coercive condition\\(see, e.g.,\cite{More}).

\begin{Th}\label{t8}
Let $K$ be a nonempty convex subset in $X$  and $T: K \rightarrow 2^{X^*}$. 
Suppose that there exist $x^0\in K$, $z^*\in T(x^0)\cap\;{\rm int}(K^*)$  and a solid convex subset $C$  in $X$ such that 

\begin{enumerate}
\item [(i)] $K(z^*,r)\subseteq {\rm int}(C)$, where $r:= \langle z^*, x^0\rangle $; 
\item [(ii)] $K \cap C \ne \emptyset $ is nonempty and compact;
\item [(iii)] $T\Big |_{K\cap C}$ is pseudomonotone;
\item [(iv)] $T\Big |_{K\cap C}$ is upper-semicontinuous, and $T(x)$ is a nonempty compact contractible subset for each $x\in K\cap C$;

Then $SOL(K, T)\cap C$ is nonempty and compact.
\end{enumerate}
\end{Th}

{\it Proof.}
Note that all assumptions of Theorem \ref{t5} are satisfied except $(iv)$. Thus, for using Theorem \ref{t5} we need to check the condition $(iv)$. Since $r = \langle z^*, x^0\rangle $ it is obvious that $x^0\in K(z^*,r)\subset K\cap {\rm int}(C)$. Suppose that $x\in K\cap \partial  C$. 
Because $K(z^*,r)\subset {\rm int}(C)$, it follows that $x \notin K(z^*,r)$. Thus, $\langle z^*, x\rangle > r = \langle z^*, x^0\rangle$. Hence, $ \langle z^*, x-x^0\rangle > 0$. Since $T\Big |_{K\cap C}$  is pseudomonotone, for each element $y^*\in T(x)$, we have $\langle y^*, x - x^0\rangle \ge 0 $, this means that the condition $(iv)$ of Theorem \ref{t5} is satisfied. Therefore, the conclusion is followed from Theorem \ref{t5}.\qed

\begin{Th}\label{t9}
Suppose that $SOL(K, T)$ is nonempty and $T$ is strictly monotone on $SOL(K, T)$. Then $SOL(K, T)$ has only one element.
\end{Th}

{\it Proof.}
Let $x_1,x_2\in SOL(K, T)$. Then, by definition,  there exist $y_1^*\in T(x_1), y_2^*\in T(x_2)$ such that  $\langle y_1^*, x - x_1\rangle \ge 0 $ and $\langle y_2^*, x - x_2\rangle \ge 0 $ for all  $x\in K$. Because $x_1,x_2\in K$, we have $\langle y_1^*, x_2 - x_1\rangle \ge 0 $ and $\langle y_2^*, x_1 - x_2\rangle \ge 0 $. Thus, $\langle y_1^*-y_2^* , x_1 - x_2\rangle = - \langle y_1^*, x_2 - x_1\rangle - \langle y_2^*, x_1 - x_2\rangle \le 0 $. Suppose that $x_1\ne x_2$. From the strict monotonicity of $T$  it follows that $\langle y_1^*-y_2^* , x_1 - x_2\rangle > 0$. The contradiction shows that $x_1= x_2$.\qed

\begin{Th}\label{t10}
Let $K$ be a nonempty convex closed subset in $X$ such that every weakly convergent sequence in $K$ strongly convergent and $T: K \rightarrow 2^{X^*}$ is upper-semiconinuous with  $T(x)$ is a nonempty compact contractible subset for each  $x\in K$. Suppose that $T$ is strongly monotone on $K$. Then $SOL(K, T)$ is a singleton.
\end{Th}

{\it Proof.}
First we prove that $SOL(K, T)$  is nonempty. If $K$ is bounded, from the closeness and convexity of $K$ in reflexive Banach space $X$, it follows that $K$ is weakly compact. Since every weakly convergent sequence in $K$ strongly convergent, we can conclude that $K$ is compact. By Lemma \ref{l4}, $SOL(K, T)$  is nonempty. If $K$ is unbounded, choose  arbitrary elements $x_0\in K$ and $y_0^*\in T(x_0)$. Because $T$ is strongly monotone, there exists a number $\alpha > 0$ such that $\langle y^*-y_0^* , x - x_0\rangle \ge \alpha ||x-x_0||^2$  for all $x\in K$  and $y^*\in T(x)$. We have
\begin{eqnarray*}
\frac {\langle y^* , x - x_0\rangle } {||x||}\ge \frac {\langle y^* , x - x_0\rangle } {||x||+ ||x_0||} & = & \frac {\langle y^*-y_0^* , x - x_0 \rangle } {||x||+ ||x_0||} + \frac {\langle y_0^* , x - x_0 \rangle} {||x||+ ||x^0||}\\
&\ge & \frac {\alpha ||x-x_0||^2} {||x||+ ||x_0||} + \frac {\langle y_0^* , x \rangle - \langle y_0^* , x_0 \rangle}{||x||+ ||x_0||}\\
&\ge & \frac {\alpha (||x||-||x_0||)^2} {||x||+ ||x_0||} + \frac {\langle y_0^* , x \rangle - \langle y_0^* , x_0 \rangle}{||x||+ ||x_0||}.\\
\end{eqnarray*}
The right side of the above inequality goes to $+ \infty$ as $||x||\rightarrow + \infty$. By Corollary \ref{c1}, $SOL(K, T)$  is nonempty. Since $T$ is strongly monotone, it is strictly monotone. By Theorem \ref{t9}, $SOL(K, T)$ is a singleton.\qed

{\bf Acknowledgments}. This work is supported by 
Vietnam's National Foundation for Science and Technology 
Development (NAFOSTED)[grant number 101.02 - 2014.51], to which, the author would
like to express many thanks. The author also thanks Professor Chuong Minh Nguyen for carefully reading and useful comments which improved the manuscript.

\end{document}